\documentclass[10pt,twoside]{article}

\usepackage{amsmath, amsthm, amssymb, mathrsfs, sansmath, setspace}

\newcommand{\ds}	{\displaystyle}

\title{Ricci Flow of Compact Locally Homogeneous Geometries on $5$-Manifolds}
\author{Thomas Bell, Amanda Hirschmann}

\begin{document}

\bibliographystyle{plain}
\maketitle

\begin{abstract}
This project serves to analyze the behavior of Ricci Flow in five dimensional manifolds. Ricci Flow was introduced by Richard Hamilton in 1982 and was an essential tool in proving the Geometrization and Poincar\'{e} Conjectures. In general, Ricci Flow is a nonlinear PDE whose solutions are rather difficult to calculate; however, in a homogeneous manifold, the Ricci Flow reduces to an ODE. The behavior of Ricci Flow in two, three, and four dimensional homogenous manifolds has been calculated and is well understood. The work presented here will describe efforts to better understand the behavior of Ricci Flow in a certain class of five dimensional homogeneous manifolds.
\end{abstract}
\section{Introduction}
Ricci Flow was first developed by Richard Hamilton in \cite{Ha} and has remained a key mathematical tool ever since. In particular, Grigori Perelman's work in \cite{Pe}, \cite{Per}, and \cite{Pere} utilized Ricci Flow in providing a proof of Thurston's Geometrization Conjecture and by extension the Poincar\'{e} Conjecture as well. Ricci Flow has been studied previously in two, three \cite{IJ}, and four \cite{IJL} dimensions as well as Backwards Ricci Flow in three \cite{CSC} and four \cite{Be} dimensions. Here we seek to analyze the same phenomenon in five dimesions; however, because of the many cases of closed five dimensional homogeneous manifolds, we restrict our analysis to only one classification, namely closed 5-dimensional contact solveable manifolds, for which the full classification can be found in \cite{DF}. These can be found in \cite{DF} which classifies all five dimensional closed homogenous manifolds. Thus, we study Ricci Flow in these manifolds by using the given equation:
\begin{equation}\label{Riccitensor}
\frac{\partial g_{ij}}{\partial t}=-2R_{ij}
\end{equation}
where $g$ is a Riemannian metric and $R$ is the associated Ricci tensor.
\\
\section{Closed five-dimensional unimodular solvable contact Lie groups}
\subsection{D1.}
Here we may choose a basis for the Lie Algebra $\{X_1,X_2,X_3,X_4,X_5\}$ such that the Lie bracket is of the form
\begin{align*}
[X_1,X_2]&=0 & [X_1,X_3]&=0 & [X_1,X_4]&=0 & [X_1,X_5]&=0 & [X_2,X_3]&=0\\
[X_2,X_4]&=X_1 & [X_2,X_5]&=0 & [X_3,X_4]&=0 & [X_3,X_5]&=X_1 & [X_4,X_5]&=0.
\end{align*}
Here we diagonalize the metric by letting $Y_i=\Lambda_i^kX_k$ with
\begin{equation*}
\Lambda=\left[\begin{array}{ccccc}
1&0&0&0&0\\a_1&1&0&0&0\\a_2&a_5&1&0&0\\a_3&a_6&a_8&1&0\\a_4&a_7&a_9&a_{10}&1
\end{array}\right].
\end{equation*}
We find that $\{Y_i\}$ satisfies the following bracket relations:
\begin{align*}
[Y_1,Y_2]&=0 & [Y_1,Y_3]&=0 & [Y_1,Y_4]&=0 & [Y_1,Y_5]&=0\\
[Y_2,Y_3]&=0 & [Y_2,Y_4]&=Y_1 & [Y_2,Y_5]&=\alpha Y_1 & [Y_3,Y_4]&=\beta Y_1\\
[Y_3,Y_5]&=(\alpha\beta+1) Y_1 & [Y_4,Y_5]&=\gamma Y_1
\end{align*}
where $\alpha=a_{10}, \beta=a_5,\text{ and }\gamma=a_6a_{10}-a_7+a_8.$
Moreover, we let
\begin{align}\label{Yhat}
\hat{Y_1}=\frac{Y_1}{\sqrt{A}}, \; \;  \hat{Y_2}=\frac{Y_2}{\sqrt{B}}, \; \; \hat{Y_3}=\frac{Y_3}{\sqrt{C}}, \; \; \hat{Y_4}=\frac{Y_4}{\sqrt{D}},\text{ and }\hat{Y_5}=\frac{Y_5}{\sqrt{E}}.
\end{align}
This gives
\begin{align*}
[\hat{Y_1},\hat{Y_2}]&=0 & [\hat{Y_1},\hat{Y_3}]&=0 & [\hat{Y_1},\hat{Y_4}]&=0\\
[\hat{Y_1},\hat{Y_5}]&=0 & [\hat{Y_2},\hat{Y_3}]&=0 & [\hat{Y_2},\hat{Y_4}]&=\frac{\sqrt{A}\hat{Y_1}}{\sqrt{BD}}\\
[\hat{Y_2},\hat{Y_5}]&=\frac{\alpha\sqrt{A}\hat{Y_1}}{\sqrt{BE}} & [\hat{Y_3},\hat{Y_4}]&=\frac{\beta\sqrt{A}\hat{Y_1}}{\sqrt{CD}} & [\hat{Y_3},\hat{Y_5}]&=\frac{(\alpha\beta+1)\sqrt{A}\hat{Y_1}}{\sqrt{CE}}\\
[\hat{Y_4},\hat{Y_5}]&=\frac{\gamma\sqrt{A}Y_1}{\sqrt{DE}}.
\end{align*}
From \cite{B}, the Ricci tensor on a Lie Group can be calculated by the equation
\begin{equation}\label{Ricci}
Ric(W,W)=-\frac{1}{2}\sum_{i}|[W,\hat{Y_i}]|^2-\frac{1}{2}\sum_{i}\langle[W,[W,\hat{Y_i}]],\hat{Y_i}\rangle+\frac{1}{2}\sum_{i<j}\langle[\hat{Y_i},\hat{Y_j}],W\rangle^2.
\end{equation}
Let $W=w_1\hat{Y_1}+w_2\hat{Y_2}+w_3\hat{Y_3}+w_4\hat{Y_4}+w_5\hat{Y_5}$. Then by (\ref{Ricci}) we get the equations of the non-zero components of the Ricci tensor:
\begin{align}\label{D1Ricci}
Ric(\hat{Y_1},\hat{Y_1})&=\frac{A}{2BD}+\frac{\alpha^2A}{2BE}+\frac{\beta^2A}{2CD}+\frac{\alpha^2\beta^2A}{2CE}+\frac{A}{2CE}+\frac{\alpha\beta A}{CE}+\frac{\gamma^2 A}{2DE}\notag\\
Ric(\hat{Y_2},\hat{Y_2})&=\frac{-A}{2BD}+\frac{-\alpha A}{2BE}\notag\\
Ric(\hat{Y_3},\hat{Y_3})&=\frac{-\beta^2}{2CD}+\frac{-\alpha^2\beta^2A}{2CE}+\frac{-A}{2CE}\notag\\
Ric(\hat{Y_4},\hat{Y_4})&=\frac{-A}{2BD}+\frac{-\beta^2A}{2CD}+\frac{-\gamma^2A}{2DE}\notag\\
Ric(\hat{Y_5},\hat{Y_5})&=\frac{-\alpha^2A}{2BE}+\frac{-\alpha\beta A}{2CE}+\frac{-A}{2CE}+\frac{-\alpha\beta A}{CE}+\frac{-\gamma^2 A}{2DE}\notag\\
Ric(\hat{Y_2},\hat{Y_3})&=\frac{-\beta A}{D(BC)^{1/2}}+\frac{-\alpha^2\beta A}{E(BC)^{1/2}}+\frac{-\alpha A}{E(BC)^{1/2}}\\
Ric(\hat{Y_2},\hat{Y_4})&=\frac{-\alpha\gamma A}{E(BD)^{1/2}}\notag\\
Ric(\hat{Y_2},\hat{Y_5})&=\frac{\gamma A}{D(BE)^{1/2}}\notag\\
Ric(\hat{Y_3},\hat{Y_4})&=\frac{-\alpha\beta\gamma A}{E(CD)^{1/2}}+\frac{-\gamma A}{E(CD)^{1/2}}\notag\\
Ric(\hat{Y_3},\hat{Y_5})&=\frac{\beta\gamma A}{D(CE)^{1/2}}\notag\\
Ric(\hat{Y_4},\hat{Y_5})&=\frac{-\alpha A}{D(CE)^{1/2}}+\frac{-\alpha\beta^2A}{C(DE)^{1/2}}+\frac{-\beta A}{C(DE)^{1/2}}.\notag
\end{align}
Then if we require $\alpha=\beta=\gamma=0$, the off-diagonal components are also forced to 0. Thus, the metric will remain diagonal under the Ricci flow and we can restate $g$.\\
\\
\indent
Let $g(t)=A(t)\theta_1^2+B(t)\theta_2^2+C(t)\theta_3^2+D(t)\theta_4^2+E(t)\theta_5^2$ where $A(0)=\lambda_1,~B(0)=\lambda_2,~C(0)=\lambda_3,~D(0)=\lambda_4,$ and $E(0)=\lambda_5$. Then Ricci flow (\ref{Riccitensor}) reduces to the following system of equations:
\begin{align}\label{D1diffeqn}
\frac{dA}{dt}&=-\frac{A^2}{BD}+-\frac{A^2}{CE}\notag\\
\frac{dB}{dt}&=\frac{A}{D}\notag\\
\frac{dC}{dt}&=\frac{A}{E}\\
\frac{dD}{dt}&=\frac{A}{B}\notag\\
\frac{dE}{dt}&=\frac{A}{C}.\notag
\end{align}
We calculate
\begin{equation}
\frac{d}{dt}(ABC)=\frac{d}{dt}(ABE)=\frac{d}{dt}(ACD)=\frac{d}{dt}(ADE)=0.
\end{equation}
Thus we have 
\begin{equation}\label{D1constants}
\begin{aligned}
ABC&=\lambda_1\lambda_2\lambda_3, & ABE&=\lambda_1\lambda_2\lambda_5,\\ ACD&=\lambda_1\lambda_3\lambda_4, & ADE&=\lambda_1\lambda_4\lambda_5.
\end{aligned}
\end{equation}
So we can restate two of the differential equations in (\ref{D1diffeqn}) as
\begin{equation}
\begin{aligned}
\frac{dB}{dt}=\frac{A}{D}=\frac{\lambda_1\lambda_2^2\lambda_3}{\lambda_4B^2C}\\
\frac{dC}{dt}=\frac{A}{E}=\frac{\lambda_1\lambda_2\lambda_3^2}{\lambda_5BC^2}.
\end{aligned}
\end{equation}
Hence,
\begin{align}\label{D1dBdC}
\frac{dB}{dC}&=\frac{\frac{dB}{dt}}{\frac{dC}{dt}}\notag\\
&=\left(\frac{\lambda_1\lambda_2^2\lambda_3}{\lambda_4B^2C}\right)\cdot\left(\frac{\lambda_5BC^2}{\lambda_1\lambda_2\lambda_3^2}\right)\notag\\
&=\frac{C}{B}\Bigg(\frac{\lambda_2\lambda_5}{\lambda_3\lambda_4}\Bigg).\\
BdB&=CdC\left(\frac{\lambda_2\lambda_5}{\lambda_3\lambda_4}\right)\notag\\
B^2&=\omega C^2+k\label{D1BC}
\end{align}
where $\ds{k=\frac{\lambda_2}{\lambda_4}\Bigg(\lambda_2\lambda_4-\lambda_3\lambda_5\Bigg)}$ and $\ds{\omega=\Bigg(\frac{\lambda_2\lambda_5}{\lambda_3\lambda_4}\Bigg)}.$
Similarly,
\begin{align}\label{D1DE}
\frac{dD}{dE}&=\frac{E}{D}\Bigg(\frac{\lambda_3\lambda_4}{\lambda_2\lambda_5}\Bigg)\notag\\
D^2&=\epsilon E^2+\ell
\end{align}
where $\ell=\ds{\frac{\lambda_4}{\lambda_2}}\Bigg(\lambda_2\lambda_4-\lambda_3\lambda_5\Bigg)$ and $\epsilon=\Bigg(\ds{\frac{\lambda_3\lambda_4}{\lambda_2\lambda_5}}\Bigg)$.
By the symmetry of $B$ and $D$ as well as $C$ and $E$ we can assume without a loss of generality that either $\lambda_2\lambda_4=\lambda_3\lambda_5$ or $\lambda_2\lambda_4<\lambda_3\lambda_5$.
\\
\\
\textbf{Case 1:} $\lambda_2\lambda_4=\lambda_3\lambda_5$
\\
\\
This gives the constants:
\begin{equation}\label{D1.1constants}
\begin{aligned}
k=\frac{\lambda_2}{\lambda_4}\left(\lambda_2\lambda_4-\lambda_3\lambda_5\right)=0\\
\ell=\frac{\lambda_4}{\lambda_2}\left(\lambda_2\lambda_4-\lambda_3\lambda_5\right)=0.
\end{aligned}
\end{equation}
Consequently,
\begin{align}\label{D1.1B}
\frac{dB}{dt}=\frac{\lambda_1\lambda_2^2\lambda_3}{\lambda_4B^2C}&=\frac{\lambda_1\lambda_2^2\lambda_3}{\lambda_4B^2}\left(\frac{1}{\omega}(B^2)\right)^{-1/2}\notag\\
B^3\frac{dB}{dt}&=\frac{\lambda_1\lambda_2^2\lambda_3\omega^{1/2}}{\lambda_4}\notag\\
B(t)&=\lambda_2\left(\frac{4\lambda_1\lambda_3\omega^{1/2}}{\lambda_2^2\lambda_4}t+1\right)^{1/4}.
\end{align}
Thus by (\ref{D1constants}) and (\ref{D1BC}) and we have the following solution to (\ref{D1diffeqn}):
\begin{align}\label{D1.1solution}
A(t)&=\lambda_1\left(\frac{4\lambda_1\lambda_3\omega^{1/2}}{\lambda_2^2\lambda_4}t+1\right)^{-1/4}\left(\frac{4\lambda_1\lambda_2}{\lambda_3^2\lambda_5\omega^{1/2}}t+1\right)^{-1/4}\notag\\
B(t)&=\lambda_2\left(\frac{4\lambda_1\lambda_3\omega^{1/2}}{\lambda_2^2\lambda_4}t+1\right)^{1/4}\notag\\
C(t)&=\lambda_3\left(\frac{4\lambda_1\lambda_2}{\lambda_3^2\lambda_5\omega^{1/2}}t+1\right)^{1/4}\\
D(t)&=\lambda_4\left(\frac{4\lambda_1\lambda_5\epsilon^{1/2}}{\lambda_2\lambda_4^2}t+1\right)^{1/4}\notag\\
E(t)&=\lambda_5\left(\frac{4\lambda_1\lambda_4}{\lambda_3\lambda_5^2\epsilon^{1/2}}t+1\right)^{1/4}.\notag
\end{align}
\\
\\
\textbf{Case 2:} $\lambda_2\lambda_4<\lambda_3\lambda_5$
\\
\\
This gives the constants:
\begin{equation}\label{D1.2constants}
\begin{aligned}
k=\frac{\lambda_2}{\lambda_4}\left(\lambda_2\lambda_4-\lambda_3\lambda_5\right)=\lambda_2^2-\frac{\lambda_2\lambda_3\lambda_5}{\lambda_4}<0\\
\ell=\frac{\lambda_4}{\lambda_2}\left(\lambda_2\lambda_4-\lambda_3\lambda_5\right)=\lambda_4^2-\frac{\lambda_3\lambda_4\lambda_5}{\lambda_2}<0.
\end{aligned}
\end{equation}
Consequently, from \cite{DF} and \cite{B} we see that
\begin{equation}\label{D1.2dB}
\begin{aligned}
\frac{dB}{dt}&=\frac{\lambda_1\lambda_2^2\lambda_3}{\lambda_4B^2}\left(\frac{1}{\omega}(B^2-k)\right)^{-1/2}\\
B^2(B^2-k)^{1/2}\frac{dB}{dt}&=\frac{\lambda_1\lambda_2^2\lambda_3\omega^{1/2}}{\lambda_4}.
\end{aligned}
\end{equation}
Which leads to,
\begin{align}
&\left(\frac{1}{8}\right)\left((B^2-k)^{1/2}(2B^3-kB)-k^2\ln\left(B+(B^2-k)^{1/2}\right)\right)\notag\\
&\qquad=\frac{\lambda_1\lambda_2^2\lambda_3\omega^{1/2}}{\lambda_4}t+k_1,
\end{align}
where $k_1$ is a constant. So for large $t$,
\begin{align}\label{D1.2B}
B(t)\approx\Bigg(\frac{4\lambda_1\lambda_2^2\lambda_3\omega^{1/2}}{\lambda_4}t\Bigg)^{1/4}.
\end{align}
Also,
\begin{equation}\label{D1.2dC}
\begin{aligned}
\frac{dC}{dt}&=\frac{\lambda_1\lambda_2\lambda_3}{\lambda_5C^2}(\omega C^2+k)^{-1/2}\\
C^2(\omega C^2+k)^{1/2}\frac{dC}{dt}&=\frac{\lambda_1\lambda_2\lambda_3^2}{\lambda_5}.
\end{aligned}
\end{equation}
Which implies,
\begin{align}
&\left(\frac{1}{8\omega^{3/2}}\right)\left((\omega^2C^2+\omega k)^{1/2}(2\omega C^3+kC)-k^2\ln(\omega C+(\omega^2C^2+\omega k)^{1/2})\right)\notag\\
&\qquad=\frac{\lambda_1\lambda_2\lambda_3^2}{\lambda_5}t+k_2,
\end{align}
where $k_2$ is a constant. So for large $t$,
\begin{align}\label{D1.2C}
C(t)\approx\Bigg(\frac{4\lambda_1\lambda_2\lambda_3^2}{\lambda_5\omega^{1/2}}t\Bigg)^{1/4}.
\end{align}
Similarly,
\begin{align}
&\left(\frac{1}{8}\right)\left((D^2-\ell)^{1/2}(2D^3-\ell D)-\ell^2\ln\left(D+(D^2-\ell)^{1/2}\right)\right)\notag\\
&\qquad=\frac{\lambda_1\lambda_4^2\lambda_5\epsilon^{1/2}}{\lambda_2}t+k_3,\\
&\left(\frac{1}{8\epsilon^{3/2}}\right)\left((\epsilon^2E^2+\epsilon\ell)^{1/2}(2\epsilon E^3+\ell E)-\ell^2\ln(\epsilon E+(\epsilon^2E^2+\epsilon\ell)^{1/2})\right)\notag\\
&\qquad=\frac{\lambda_1\lambda_4\lambda_5^2}{\lambda_3}t+k_4,
\end{align}
where $k_3$ and $k_4$ are constants. So for large $t$,
\begin{align}
D(t)&\approx\Bigg(\frac{4\lambda_1\lambda_4^2\lambda_5\epsilon^{1/2}}{\lambda_2}t\Bigg)^{1/4}\label{D1.2D}\\
E(t)&\approx\Bigg(\frac{4\lambda_1\lambda_4\lambda_5^2}{\lambda_3\epsilon^{1/2}}t\Bigg)^{1/4}.\label{D1.2E}
\end{align}
Then by (\ref{D1constants}) we have the following solution to (\ref{D1diffeqn}) for large values of $t$:
\begin{align}
A(t)&\approx\frac{1}{2t^{1/2}}\Bigg(\lambda_1^2\lambda_2\lambda_3\lambda_4\lambda_5\Bigg)^{1/4}\notag\\
B(t)&\approx\Bigg(\frac{4\lambda_1\lambda_2^2\lambda_3\omega^{1/2}}{\lambda_4}t\Bigg)^{1/4}\notag\\
C(t)&\approx\Bigg(\frac{4\lambda_1\lambda_2\lambda_3^2}{\lambda_5\omega^{1/2}}t\Bigg)^{1/4}\\
D(t)&\approx\Bigg(\frac{4\lambda_1\lambda_4^2\lambda_5\epsilon^{1/2}}{\lambda_2}t\Bigg)^{1/4}\notag\\
E(t)&\approx\Bigg(\frac{4\lambda_1\lambda_4\lambda_5^2}{\lambda_3\epsilon^{1/2}}t\Bigg)^{1/4}\notag.
\end{align}

\subsection{D2.}

Here we may choose a basis for the Lie Algebra $\{X_1,X_2,X_3,X_4,X_5\}$ such that the Lie bracket is of the form
\begin{align*}
[X_1,X_2]&=0 & [X_1,X_3]&=0 & [X_1,X_4]&=0 & [X_1,X_5]&=0 & [X_2,X_3]&=0\\
[X_2,X_4]&=0 & [X_2,X_5]&=X_1 & [X_3,X_4]&=X_1 & [X_3,X_5]&=X_2 & [X_4,X_5]&=0.
\end{align*}
\\
Here we diagonalize the metric by letting $Y_i=\Lambda_i^kX_k$ with
\begin{equation*}
\Lambda=\left[\begin{array}{ccccc}
1&0&0&0&0\\a_1&1&0&0&0\\a_2&a_5&1&0&0\\a_3&a_6&a_8&1&0\\a_4&a_7&a_9&a_{10}&1
\end{array}\right].
\end{equation*}
\\
We find that $\{Y_i\}$ satisfies the following bracket relations:
\begin{align*}
[Y_1,Y_2]&=0 & [Y_1,Y_3]&=0 & [Y_1,Y_4]&=0 & [Y_1,Y_5]&=0\\
[Y_2,Y_3]&=0 & [Y_2,Y_4]&=0 & [Y_2,Y_5]&=Y_1 & [Y_3,Y_4]&=Y_1\\
[Y_3,Y_5]&=\alpha Y_1+Y_2 & [Y_4,Y_5]&=\beta Y_1+\gamma Y_2
\end{align*}
where $\alpha=a_{10}+a_5-a_1, \beta=a_8a_{10}-a_9+a_6-a_1a_8,\text{ and }\gamma=a_8.$\\
\\
\indent
Define $\ds{\hat{Y_1}=\frac{Y_1}{\sqrt{A}}}$ as in (\ref{Yhat}). Let $W=w_1\hat{Y_1}+w_2\hat{Y_2}+w_3\hat{Y_3}+w_4\hat{Y_4}+w_5\hat{Y_5}$. Then by (\ref{Ricci}) we get the equations of the non-zero components of the Ricci tensor:
\begin{equation}
\begin{aligned}\label{D2Ricci}
Ric(\hat{Y_1},\hat{Y_1})&=\frac{A}{2BE}+\frac{A}{2CD}+\frac{\alpha^2A}{2CE}+\frac{\beta^2A}{2DE}\\
Ric(\hat{Y_2},\hat{Y_2})&=\frac{B}{2CE}+\frac{\gamma^2A}{2DE}+\frac{-A}{2BE}\\
Ric(\hat{Y_3},\hat{Y_3})&=\frac{-A}{2CD}+\frac{2\alpha^2A}{2CE}+\frac{-B}{2CE}\\
Ric(\hat{Y_4},\hat{Y_4})&=\frac{-A}{2CD}+\frac{-\beta^2A}{2DE}+\frac{-\gamma^2A}{2DE}\\
Ric(\hat{Y_5},\hat{Y_5})&=\frac{-A}{2BE}+\frac{-\alpha^2A}{2CE}+\frac{-B}{2CE}+\frac{-\beta^2A}{2DE}+\frac{-\gamma^2B}{2DE}\\
Ric(\hat{Y_2},\hat{Y_3})&=\frac{-\alpha A}{E(BC)^{1/2}}\\
Ric(\hat{Y_2},\hat{Y_4})&=\frac{-\beta A}{E(BD)^{1/2}}\\
Ric(\hat{Y_3},\hat{Y_4})&=\frac{-\alpha\beta A-\gamma B}{E(CD)^{1/2}}\\
Ric(\hat{Y_3},\hat{Y_5})&=\frac{\beta A}{D(CE)^{1/2}}\\
Ric(\hat{Y_4},\hat{Y_5})&=\frac{-\alpha A}{C(DE)^{1/2}}.
\end{aligned}
\end{equation}
Then if we require $\alpha=\beta=\gamma=0$, the off-diagonal components are also forced to 0. Thus, the metric will remain diagonal under the Ricci flow and we can restate $g$.\\
\\
\indent
Let $g(t)=A(t)\theta_1^2+B(t)\theta_2^2+C(t)\theta_3^2+D(t)\theta_4^2+E(t)\theta_5^2$ where $A(0)=\lambda_1,~B(0)=\lambda_2,~C(0)=\lambda_3,~D(0)=\lambda_4,$ and $E(0)=\lambda_5$. Then Ricci flow (\ref{Riccitensor}) reduces to the following system of equations:
\begin{align}\label{D2diffeqn}
\frac{dA}{dt}&=-\frac{A^2}{BE}+-\frac{A^2}{CD}\notag\\
\frac{dB}{dt}&=-\frac{B^2}{CE}+\frac{A}{E}\notag\\
\frac{dC}{dt}&=\frac{A}{D}+\frac{B}{E}\\
\frac{dD}{dt}&=\frac{A}{C}\notag\\
\frac{dE}{dt}&=\frac{A}{B}+\frac{B}{C}.\notag
\end{align}
We calculate
\begin{equation}
\frac{d}{dt}(ABC)=\frac{d}{dt}(A^2BD^2E)=0.
\end{equation}
Thus we have
\begin{equation}\label{D2constants}
\begin{aligned}
ABC&=\lambda_1\lambda_2\lambda_3\\
A^2BD^2E&=\lambda_1^2\lambda_2\lambda_4^2\lambda_5.
\end{aligned}
\end{equation}
So we can restate three of the differential equations in (\ref{D2diffeqn}) as
\begin{equation}\label{D2diffeqn1}
\begin{aligned}
\frac{dA}{dt}&=\frac{-A^3}{\lambda_1\lambda_2}\left(\frac{AD^2}{\lambda_1\lambda_4^2\lambda_5}+\frac{B}{\lambda_3D}\right)\\
\frac{dB}{dt}&=\frac{A^3D^2}{\lambda_1^2\lambda_2\lambda_4^2\lambda_5}\left(\frac{-B^4}{\lambda_1\lambda_2\lambda_3}+B\right)\\
\frac{dD}{dt}&=\frac{A^2B}{\lambda_1\lambda_2\lambda_3}.
\end{aligned}
\end{equation}
Then we calculate
\begin{align}\label{D2ACB^2}
\frac{d}{dt}\Bigg(\frac{AC}{B^2}\Bigg)&=\frac{3A}{B^2E}(B^2-AC)\notag\\
&=3\Bigg(\frac{AC}{B^2}\Bigg)\Bigg(\frac{B^2-AC}{BCE}\Bigg).
\end{align}
Hence, we have demonstrated that $\ds{\frac{AC}{B^2}\rightarrow1}$ and we can consider two cases: $AC=B^2$ or $AC\not=B^2$.
\\
\\
\textbf{Case 1:} $AC=B^2$
\\
\\
If $\ds{B^2=AC}$, then $\ds{\frac{dB}{dt}=0}$ and $B$ is a constant. In fact, by ($\ref{D2constants}$) we figure
\begin{align}
B^2&=AC\notag\\
B^3&=\lambda_1\lambda_2\lambda_3\\
B&=\lambda_2=(\lambda_1\lambda_2\lambda_3)^{1/3}.\label{D2constantB}
\end{align}
Then follows,
\begin{align}\label{D2dAdD}
\frac{dA}{dD}&=\frac{\frac{dA}{dt}}{\frac{dD}{dt}}\notag\\
&=\left(\frac{-A^4D^2}{\lambda_1^2\lambda_2\lambda_4^2\lambda_5}+\frac{-A^3B}{\lambda_1\lambda_2\lambda_3D}\right)\left(\frac{\lambda_1\lambda_2\lambda_3}{A^2B}\right)\notag\\
&=\frac{-\lambda_3A^2D^2}{\lambda_1\lambda_4^2\lambda_5B}+\frac{-A}{D}\notag\\
&=\frac{-\lambda_3A^2D^2}{\lambda_1\lambda_4^2\lambda_5(\lambda_1\lambda_2\lambda_3)^{1/3}}+\frac{-A}{D}.
\end{align}
Let $\ds{\ell=\frac{\lambda_3}{\lambda_1\lambda_4^2\lambda_5(\lambda_1\lambda_2\lambda_3)^{1/3}}}$, then we can restate (\ref{D2dAdD}) as
\begin{equation}
\begin{aligned}
\frac{dA}{dD}+\frac{A}{D}=-\ell A^2D^2.
\end{aligned}
\end{equation}
This gives us a Bernoulli Equation which we can calculate as in \cite[p.77]{BD}. We can solve by replacing with $v=A^{-1}$:
\begin{equation}\label{D2dvdD}
\begin{aligned}
\frac{dv}{dD}&+(-D^{-1})v=\ell D^2.
\end{aligned}
\end{equation}
We can find the solution to this first-order linear differential equation by using an integrating factor:
\begin{align}\label{D2intfactor}
v&=\left(\frac{(-1)\int e^{-\int D^{-1}dD}(-\ell D^2)dD}{e^{-\int D^{-1}dD}}\right)\notag\\
&=\frac{\ell}{2}D^3+KD\notag\\
&=\frac{\ell}{2}(D^3+K_1D)
\end{align}
where $K$ is a constant and $\ds{K_1=\frac{2K}{\ell}}$. Because $A=v^{-1}$, we know that
\begin{align}
A&=\frac{2}{\ell}(D^3+K_1D)^{-1}\notag\\
&=m(D^3+K_1D)^{-1}
\end{align}
where $\ds{m=\frac{2}{\ell}}$.
Furthermore, by (\ref{D2diffeqn1})
\begin{align}
\frac{dD}{dt}&=\frac{m^2B}{\lambda_1\lambda_2\lambda_3}(D^3+K_1D)^{-2}\notag\\
&=\frac{m^2(\lambda_1\lambda_2\lambda_3)^{1/3}}{\lambda_1\lambda_2\lambda_3}(D^3+K_1D)^{-2}.
\intertext{Let $\ds{M=\frac{m^2}{(\lambda_1\lambda_2\lambda_3)^{2/3}}}$, and then}
\frac{dD}{dt}&=M(D^3+K_1D)^{-2}\notag\\
(D^3+K_1D)^2\frac{dD}{dt}&=M.\notag
\end{align}
With this we can see that for large $t$,
\begin{align}\label{D2D}
D^7&\approx(Mt)\notag\\
D&\approx(Mt)^{1/7}.
\intertext{Furthermore,}
A&=m(D^3+K_1D)^{-1}\notag\label{D2A}\\
&\approx\frac{m}{D^3}\notag\\
&\approx m(Mt)^{-3/7}.
\end{align}
We use this result and (\ref{D2constants}) to calculate
\begin{align}\label{D2C}
C&=\frac{\lambda_1\lambda_2\lambda_3}{AB}\notag\\
&\approx\frac{(\lambda_1\lambda_2\lambda_3)^{2/3}}{m}(Mt)^{3/7}\notag\\
&=M_C(Mt)^{3/7}.
\end{align}
The results in (\ref{D2constantB}) and (\ref{D2A}) as well as (\ref{D2constants}) enable us to solve for $E$:
\begin{align}
E&=\frac{\lambda_1^2\lambda_2\lambda_4^2\lambda_5}{A^2BD^2}\notag\\
&\approx\frac{\lambda_1^2\lambda_2\lambda_4^2\lambda_5(Mt)^{6/7}}{m^2(\lambda_1\lambda_2\lambda_3)^{1/3}(Mt)^{2/7}}\notag\\
&=M_E(Mt)^{4/7}.\label{D2E}
\end{align}
Then from (\ref{D2A}), (\ref{D2constantB}), (\ref{D2C}), (\ref{D2D}), and (\ref{D2E}) we have the following solution to (\ref{D2diffeqn}) for large values of $t$:
\begin{align}\label{D2.1solution}
A(t)&\approx m(Mt)^{-3/7}\notag\\
B(t)&=(\lambda_1\lambda_2\lambda_3)^{1/3}\notag\\
C(t)&\approx M_C(Mt)^{3/7}\\
D(t)&\approx(Mt)^{1/7}\notag\\
E(t)&\approx M_E(Mt)^{4/7}.\notag
\end{align}
\\
\\
\textbf{Case 2:} $AC\not=B^2$
\\
\\
We know that $\ds{\frac{AC}{B^2}\rightarrow1}$. We analyze the behavior of B, specifically we want to know the value of
\begin{equation}\label{D2.2B}
\frac{-B^4}{\lambda_1\lambda_2\lambda_3}+B.
\end{equation}
Therefore, we calculate said behavior.
\begin{equation}\label{D2.2constants}
\begin{aligned}
\text{If }\ds{\frac{-B^4}{\lambda_1\lambda_2\lambda_3}+B>0},\text{ then }\ds{\frac{dB}{dt}>0}\text{ and $B$ is increasing.}\\
\text{If }\ds{\frac{-B^4}{\lambda_1\lambda_2\lambda_3}+B<0},\text{ then }\ds{\frac{dB}{dt}<0}\text{ and $B$ is decreasing.}\\
\text{If }\ds{\frac{-B^4}{\lambda_1\lambda_2\lambda_3}+B=0},\text{ then }\ds{\frac{dB}{dt}=0}\text{ and $B$ is a constant.}
\end{aligned}
\end{equation}
$\ds{\frac{-B^4}{\lambda_1\lambda_2\lambda_3}+B>0}$ and $\ds{\frac{-B^4}{\lambda_1\lambda_2\lambda_3}+B<0}$ indicate that $B>(\lambda_1\lambda_2\lambda_3)^{1/3}$ and $B<(\lambda_1\lambda_2\lambda_3)^{1/3}$, respectively. Hence, we may conclude that $B$ is approaching the constant $(\lambda_1\lambda_2\lambda_3)^{1/3}$. We say that
\begin{equation}\label{D2.2constantB}
0<b_1\leq B\leq b_2.
\end{equation}
Then follows,
\begin{align}\label{D2.2dAdD}
\frac{dA}{dD}&=\left(\frac{-A^4D^2}{\lambda_1^2\lambda_2\lambda_4^2\lambda_5}+\frac{-A^3B}{\lambda_1\lambda_2\lambda_3D}\right)\left(\frac{\lambda_1\lambda_2\lambda_3}{A^2B}\right)\notag\\
&=\frac{-\lambda_3A^2D^2}{\lambda_1\lambda_4^2\lambda_5B}+\frac{-A}{D}
\end{align}
\begin{equation}\label{D2.2dAdDinequality}
\frac{-\lambda_3A^2D^2}{\lambda_1\lambda_4^2\lambda_5b_2}+\frac{-A}{D}\leq\frac{dA}{dD}\leq\frac{-\lambda_3A^2D^2}{\lambda_1\lambda_4^2\lambda_5b_1}+\frac{-A}{D}.
\end{equation}
Let $\ds{\frac{\lambda_3}{\lambda_1\lambda_4^2\lambda_5b_1}=\ell_1}$ and $\ds{\frac{\lambda_3}{\lambda_1\lambda_4^2\lambda_5b_2}=\ell_2},$ then
\begin{align}
\frac{dA}{dD}+\frac{A}{D}\leq-\ell_1A^2D^2.
\end{align}
We have now developed two first-order linear differential equations and we can solve them by proceeding as we did in Case 1. Then by (\ref{D2constants}) and since $B\rightarrow(\lambda_1\lambda_2\lambda_3)^{1/3}$ for large values of $t$ we have a nearly identical result to (\ref{D2.1solution}) as found in Case 1:
\begin{align}\label{D2.2solution}
A(t)&\approx M_1t^{-3/7}\notag\\
B(t)&\rightarrow(\lambda_1\lambda_2\lambda_3)^{1/3}\notag\\
C(t)&\approx M_3t^{3/7}\\
D(t)&\approx M_4t^{1/7}\notag\\
E(t)&\approx M_5t^{4/7}.\notag
\end{align}
\subsection{D3.}

Here we may choose a basis for the Lie Algebra $\{X_1,X_2,X_3,X_4,X_5\}$ such that the Lie bracket is of the form
\begin{align*}
[X_1,X_2]&=0 & [X_1,X_3]&=0 & [X_1,X_4]&=0 & [X_1,X_5]&=0 & [X_2,X_3]&=0\\
[X_2,X_4]&=0 & [X_2,X_5]&=X_1 & [X_3,X_4]&=X_1 & [X_3,X_5]&=X_2 & [X_4,X_5]&=X_3.
\end{align*}
Here we diagonalize the metric by letting $Y_i=\Lambda_i^kX_k$ with
\begin{equation*}
\Lambda=\left[\begin{array}{ccccc}
1&0&0&0&0\\a_1&1&0&0&0\\a_2&a_5&1&0&0\\a_3&a_6&a_8&1&0\\a_4&a_7&a_9&a_{10}&1
\end{array}\right].
\end{equation*}
We find that $\{Y_i\}$ satisfies the following bracket relations:
\begin{align*}
[Y_1,Y_2]&=0 & [Y_1,Y_3]&=0 & [Y_1,Y_4]&=0 & [Y_1,Y_5]&=0\\
[Y_2,Y_3]&=0 & [Y_2,Y_4]&=0 & [Y_2,Y_5]&=Y_1 & [Y_3,Y_4]&=Y_1\\
[Y_3,Y_5]&=\alpha Y_1+Y_2-Y_3 & [Y_4,Y_5]&=\beta Y_1+\gamma Y_2+Y_3
\end{align*}
where $\alpha=a_{10}+a_5-a_1, \beta=a_8a_{10}-a_9+a_6-a_1a_8-a_2+a_1a_5,\text{ and }\gamma=-a_5+a_{10}.$\\
\\
\indent
Define $\ds{\hat{Y_1}=\frac{Y_1}{\sqrt{A}}}$ as in (\ref{Yhat}). Let $W=w_1\hat{Y_1}+w_2\hat{Y_2}+w_3\hat{Y_3}+w_4\hat{Y_4}+w_5\hat{Y_5}$. Then by (\ref{Ricci}) we get the equations of the non-zero components of the Ricci tensor:
\begin{equation}\label{D3Ricci}
\begin{aligned}
Ric(\hat{Y_1},\hat{Y_1})&=\frac{A}{2BE}+\frac{A}{2CD}+\frac{\alpha^2A}{2CE}+\frac{\beta^2A}{2DE}\\
Ric(\hat{Y_2},\hat{Y_2})&=\frac{B}{2CE}+\frac{\gamma^2B}{2DE}+\frac{-A}{2BE}\\
Ric(\hat{Y_3},\hat{Y_3})&=\frac{C}{2DE}+\frac{-A}{2CD}+\frac{-\alpha^2A}{2CE}+\frac{-B}{2CE}\\
Ric(\hat{Y_4},\hat{Y_4})&=\frac{-A}{2CD}+\frac{-\beta^2A}{2DE}+\frac{-\gamma^2B}{2DE}+\frac{-C}{2DE}\\
Ric(\hat{Y_5},\hat{Y_5})&=\frac{-A}{2BE}+\frac{-\alpha^2A}{2CE}+\frac{-B}{2CE}+\frac{-\beta^2A}{2DE}+\frac{-\gamma^2B}{2DE}+\frac{-C}{2DE}\\
Ric(\hat{Y_2},\hat{Y_3})&=\frac{-\alpha A}{E(BC)^{1/2}}\\
Ric(\hat{Y_2},\hat{Y_4})&=\frac{-\beta A}{E(BD)^{1/2}}\\
Ric(\hat{Y_3},\hat{Y_4})&=\frac{-\alpha\beta A-\gamma B}{E(CD)^{1/2}}\\
Ric(\hat{Y_3},\hat{Y_5})&=\frac{-\beta A}{D(CE)^{1/2}}\\
Ric(\hat{Y_4},\hat{Y_5})&=\frac{-\alpha A}{C(DE)^{1/2}}.
\end{aligned}
\end{equation}
Then if we require $\alpha=\beta=\gamma=0$, the off-diagonal components are also forced to 0. Thus, the metric will remain diagonal under the Ricci flow and we can restate $g$.\\
\\
\indent
Let $g(t)=A(t)\theta_1^2+B(t)\theta_2^2+C(t)\theta_3^2+D(t)\theta_4^2+E(t)\theta_5^2$ where $A(0)=\lambda_1,~B(0)=\lambda_2,~C(0)=\lambda_3,~D(0)=\lambda_4,$ and $E(0)=\lambda_5$. Then Ricci flow (\ref{Riccitensor}) reduces to the following system of equations:
\begin{align}\label{D3diffeqn}
\frac{dA}{dt}&=-\frac{A^2}{BE}+-\frac{A^2}{CD}\notag\\
\frac{dB}{dt}&=-\frac{B^2}{CE}+\frac{A}{E}\notag\\
\frac{dC}{dt}&=-\frac{C^2}{DE}+\frac{A}{D}+\frac{B}{E}\\
\frac{dD}{dt}&=\frac{A}{C}+\frac{C}{E}\notag\\
\frac{dE}{dt}&=\frac{A}{B}+\frac{B}{C}+\frac{C}{D}.\notag
\end{align}
We calculate
\begin{equation*}\label{D3constants}
\frac{d}{dt}(A^5B^4C^3D^2E)=0.
\end{equation*}
Moreover, we calculate a new system of differential equations
\begin{equation}\label{D3diffeqn1}
\begin{aligned}
\frac{d}{dt}\left(\frac{A}{BE}\right)&=-3\left(\frac{A}{BE}\right)^2+\left(\frac{A}{BE}\right)\left(-\frac{A}{CD}\right)+\left(\frac{A}{BE}\right)\left(-\frac{C}{DE}\right)\\
\frac{d}{dt}\left(\frac{A}{CD}\right)&=-3\left(\frac{A}{CD}\right)^2+\left(\frac{A}{CD}\right)\left(-\frac{A}{BE}\right)+\left(\frac{A}{CD}\right)\left(-\frac{B}{CE}\right)\\
\frac{d}{dt}\left(\frac{B}{CE}\right)&=-3\left(\frac{B}{CE}\right)^2+\left(\frac{B}{CE}\right)\left(-\frac{A}{CD}\right)\\
\frac{d}{dt}\left(\frac{C}{DE}\right)&=-3\left(\frac{C}{DE}\right)^2+\left(\frac{C}{DE}\right)\left(-\frac{A}{BE}\right).
\end{aligned}
\end{equation}

Then we let
\begin{equation}\label{D2constants1}
x=\frac{A}{BE}, y=\frac{A}{CD},z=\frac{B}{CE},\text{ and }w=\frac{C}{DE},
\end{equation}

and (\ref{D3diffeqn1}) becomes
\begin{equation}\label{D3diffeqn2}
\begin{aligned}
x'&=-3x^2-xy-xw\\
y'&=-3y^2-xy-yz\\
z'&=-3z^2-yz\\
w'&=-3w^2-xw.
\end{aligned}
\end{equation}
Consequently, we may restate (\ref{D3diffeqn}) as
\begin{equation}\label{D3diffeqn3}
\begin{aligned}
\frac{1}{A}\frac{dA}{dt}&=-x-y\\
\frac{1}{B}\frac{dB}{dt}&=-z+x\\
\frac{1}{C}\frac{dC}{dt}&=-w+y+z\\
\frac{1}{D}\frac{dD}{dt}&=y+w\\
\frac{1}{E}\frac{dE}{dt}&=x+z+w.
\end{aligned}
\end{equation}
Moreover, we know by (\ref{D3diffeqn2}) that $x'=-3x^2-xy-xw$.
Assume for the moment that $\ds{x'=-kx^2}$ where  $k$ is a constant. Hence, $\ds{x=\frac{k_1}{(t+c)}}$ where $c$ is a constant and $\ds{k_1=\frac{1}{k}}$. This leads to the further assumptions:
\begin{align}
&x=\frac{k_1}{t+c}, y=\frac{k_2}{t+c}, z=\frac{k_3}{t+c},\text{ and }w=\frac{k_4}{t+c}\label{D3constantsk}
\intertext{which, in turn, imply}
&x'=\frac{-k_1}{(t+c)^2}, y'=\frac{-k_2}{(t+c)^2}, z'=\frac{-k_3}{(t+c)^2},\text{ and }w'=\frac{-k_4}{(t+c)^2}\label{D3derivativesk}
\end{align}
with the caveat that $c$ has the same value in $x, y, z,$ and $w$.
If we substitute these equations into (\ref{D3diffeqn2}) we  get the following set of equations:
\begin{equation}\label{D3diffeqn4}
\begin{aligned}
1&=3k_1+k_2+k_4\\
1&=3k_2+k_1+k_3\\
1&=3k_3+k_2\\
1&=3k_4+k_1.
\end{aligned}
\end{equation}
This gives the solutions
\begin{equation}\label{D3constants3}
k_1=\frac{2}{11}, k_2=\frac{2}{11}, k_3=\frac{3}{11},\text{ and }k_4=\frac{3}{11}.
\end{equation}
Therefore, by (\ref{D3diffeqn3}),
\begin{align}\label{D3A}
\frac{1}{A}\frac{dA}{dt}&=-\frac{A}{BE}-\frac{A}{CD}\notag\\
&=\frac{-\frac{4}{11}}{t+c}\notag\\
\ln A&=-\frac{4}{11}\ln(t+c)+\ell_1\notag\\
A&=\ell_A(t+c)^{-4/11}
\end{align}
where $\ell_1$ and $\ell_A$ are constants. Similarly,
\begin{equation}\label{D3solution}
\begin{aligned}
A&=\ell_A(t+c)^{-4/11}\\
B&=\ell_B(t+c)^{-1/11}\\
C&=\ell_C(t+c)^{2/11}\\
D&=\ell_D(t+c)^{5/11}\\
E&=\ell_E(t+c)^{8/11}.
\end{aligned}
\end{equation}
By (\ref{D3constants3}) we see that
\begin{align}
3x=3y&=2z=2w\label{D3initialconditions}\\
3\bigg(\frac{A}{BE}\bigg)=3\bigg(\frac{A}{CD}\bigg)&=2\bigg(\frac{B}{CE}\bigg)=2\bigg(\frac{C}{DE}\bigg).\label{D3initialconditions1}\\
\intertext{Plugging in $t=0$ gives the initial requirements for this set of solutions:}\notag
\lambda_2\lambda_5=\lambda_3\lambda_4, \lambda_2\lambda_4&=\lambda_3^2,\text{ and }3\lambda_1\lambda_5=2\lambda_3^2.
\end{align}
Now we calculate (\ref{D3diffeqn2}) under more general circumstances. First we solve:
\begin{align}
\bigg(\frac{x}{y}\bigg)'&=\bigg(\frac{x}{y}\bigg)\Big((2y-2x)+(z-w)\Big)\\
\bigg(\frac{z}{w}\bigg)'&=\bigg(\frac{z}{w}\bigg)\Big((3w-3z)+(y-x)\Big)\\
\bigg(\frac{x}{z}\bigg)'&=\bigg(\frac{x}{z}\bigg)\Big((2z-3x)+(z-w)\Big)\\
\bigg(\frac{x}{w}\bigg)'&=\bigg(\frac{x}{w}\bigg)\Big((2w-3x)+(x-y)\Big)\\
\bigg(\frac{y}{z}\bigg)'&=\bigg(\frac{y}{z}\bigg)\Big((2z-3y)+(y-x)\Big)\\
\bigg(\frac{y}{w}\bigg)'&=\bigg(\frac{y}{w}\bigg)\Big((2w-3y)+(w-z)\Big)\\
\bigg(\frac{x}{y}\bigg)'&=\bigg(\frac{x}{y}\bigg)\Big((2y-2x)+(z-w)\Big)\\
\bigg(\frac{z}{w}\bigg)'&=\bigg(\frac{z}{w}\bigg)\Big((3w-3z)+(y-x)\Big)
\end{align}
If $\ds{\bigg(\frac{x}{y}\bigg)}$ is moving away from $1$, then without a loss of generality let $x\geq y$. Thus $\ds{(2y-2x)\leq0}$, but $\ds{z-w>2(x-y)\geq0}$, so $\ds{\frac{z}{w}>1}$ and $\ds{3(w-z)+(y-x)<\frac{7}{2}(w-z)}$. So proportionally, $\ds{\bigg(\frac{z}{w}\bigg)}$ is approaching $1$ by at least $\ds{\bigg(\frac{7}{2}\bigg)}$ the rate that $\ds{\bigg(\frac{x}{y}\bigg)}$ is moving away from $1$. Likewise, if $\ds{\bigg(\frac{z}{w}\bigg)}$ is moving away from $1$, then $\ds{\bigg(\frac{x}{y}\bigg)}$ is moving towards $1$, proportionally, by at least $\ds{\bigg(\frac{5}{2}\bigg)}$ the rate that $\ds{\bigg(\frac{z}{w}\bigg)}$ is moving away. Since $x,y,w,$ and $z$ all approach $0$, then $\ds{\bigg(\frac{x}{y}\bigg)}$ and $\ds{\bigg(\frac{z}{w}\bigg)}$ must approach $1$. Now we see that $\ds{\bigg(\frac{x}{z}\bigg)}$, $\ds{\bigg(\frac{x}{w}\bigg)}$, $\ds{\bigg(\frac{y}{z}\bigg)}$, and $\ds{\bigg(\frac{y}{w}\bigg)}$ must all approach $\ds{\frac{2}{3}}$. Thus, by (\ref{D3diffeqn2}) we see that for large $t$,
\begin{equation}\label{D3diffeqn5}
\begin{aligned}
x'&\approx-\frac{11}{2}x^2, & y'&\approx-\frac{11}{2}y^2,\\
z'&\approx-\frac{11}{3}z^2, & w'&\approx-\frac{11}{3}w^2.
\end{aligned}
\end{equation}
Therefore, we have the following solution to (\ref{D3diffeqn}),
\begin{align}
A&\approx\ell_At^{-4/11}\notag\\
B&\approx\ell_Bt^{-1/11}\notag\\
C&\approx\ell_Ct^{2/11}\\
D&\approx\ell_Dt^{5/11}\notag\\
E&\approx\ell_Et^{8/11}.\notag
\end{align}
\subsection{D5.}
Here we may choose a basis for the Lie Algebra $\{X_1,X_2,X_3,X_4,X_5\}$ such that the Lie bracket is of the form
\begin{align*}
[X_1,X_2]&=0 & [X_1,X_3]&=0 & [X_1,X_4]&=0 & [X_1,X_5]&=0 & [X_2,X_3]&=X_1\\
[X_2,X_4]&=0 & [X_2,X_5]&=X_2 & [X_3,X_4]&=0 & [X_3,X_5]&=-X_3 & [X_4,X_5]&=X_1.
\end{align*}
Here we diagonalize the metric by letting $Y_i=\Lambda_i^kX_k$ with
\begin{equation*}
\Lambda=\left[\begin{array}{ccccc}
1&0&0&0&0\\a_1&1&0&0&0\\a_2&a_5&1&0&0\\a_3&a_6&a_8&1&0\\a_4&a_7&a_9&a_{10}&1
\end{array}\right].
\end{equation*}
We find that $\{Y_i\}$ satisfies the following bracket relations:
\begin{align*}
[Y_1,Y_2]&=0 & [Y_1,Y_3]&=0 & [Y_1,Y_4]&=0\\
[Y_1,Y_5]&=0 & [Y_2,Y_3]&=Y_1 & [Y_2,Y_4]&=\alpha Y_1\\
[Y_2,Y_5]&=\beta Y_1+Y_2 & [Y_3,Y_4]&=\gamma Y_1 & [Y_3,Y_5]&=\delta Y_1+\eta Y_2-Y_3\\
[Y_4,Y_5]&=\mu Y_1+\rho Y_2-\alpha Y_3
\end{align*}
where $\alpha=a_{10}, \beta=a_6-a_1, \gamma=-a_5a_{10}-a_6, \delta=a_5a_6-a_7+a_2-2a_1a_5, \eta=2a_5, \mu=a_6a_9-a_7a_8+1-a_1a_6+a_2a_{10}-a_1a_5a_{10},\text{ and }\rho=a_6+a_5a_{10}.$\\
\\
\indent
Define $\ds{\hat{Y_1}=\frac{Y_1}{\sqrt{A}}}$ as in (\ref{Yhat}). Let $W=w_1\hat{Y_1}+w_2\hat{Y_2}+w_3\hat{Y_3}+w_4\hat{Y_4}+w_5\hat{Y_5}$. Then by (\ref{Ricci}) we get the equations of the non-zero components of the Ricci tensor:
\begin{equation}\label{D5Ricci}
\begin{aligned}
Ric(\hat{Y_1},\hat{Y_1})&=\frac{A}{2BC}+\frac{\alpha^2A}{2BD}+\frac{\beta^2A}{2BE}+\frac{\gamma^2A}{2CD}+\frac{\delta^2A}{2CE}+\frac{\mu^2A}{2DE}\\
Ric(\hat{Y_2},\hat{Y_2})&=\frac{\eta^2B}{2CE}+\frac{\rho^2B}{2DE}+\frac{-A}{2BC}+\frac{-\alpha^2A}{2BD}+\frac{-\beta^2A}{2BE}\\
Ric(\hat{Y_3},\hat{Y_3})&=\frac{\alpha^2C}{2DE}+\frac{-A}{2BC}+\frac{-\gamma^2A}{2CD}+\frac{-\delta^2A}{2CE}+\frac{-\eta^2B}{2CE}\\
Ric(\hat{Y_4},\hat{Y_4})&=\frac{-\alpha^2A}{2BD}+\frac{-\gamma^2A}{2CD}+\frac{-\mu^2A}{2DE}+\frac{-\rho^2B}{2DE}+\frac{-\alpha^2C}{2DE}\\
Ric(\hat{Y_5},\hat{Y_5})&=\frac{-\beta^2A}{2DE}+\frac{-\delta^2A}{2CE}+\frac{-\eta^2B}{2CE}+\frac{-\mu^2A}{2DE}+\frac{-\rho^2B}{2DE}+\frac{-\alpha^2C}{2DE}\\
&+\frac{-2}{E}\\
Ric(\hat{Y_2},\hat{Y_3})&=\frac{-\alpha\gamma A}{D(BC)^{1/2}}+\frac{-\beta\delta A}{E(BC)^{1/2}}+\frac{-\eta B}{E(BC)^{1/2}}\\
Ric(\hat{Y_2},\hat{Y_4})&=\frac{\gamma A}{C(BD)^{1/2}}+\frac{-\beta\mu A}{E(BD)^{1/2}}+\frac{-\rho B}{E(BD)^{1/2}}\\
Ric(\hat{Y_2},\hat{Y_5})&=\frac{\delta A}{C(BE)^{1/2}}+\frac{\alpha\mu A}{D(BE)^{1/2}}\\
Ric(\hat{Y_3},\hat{Y_4})&=\frac{-\alpha A}{B(CD)^{1/2}}+\frac{-\delta\mu A}{E(CD)^{1/2}}+\frac{-\eta\rho B}{E(CD)^{1/2}}+\frac{-\alpha C}{E(CD)^{1/2}}\\
Ric(\hat{Y_3},\hat{Y_5})&=\frac{-\beta A}{B(CE)^{1/2}}+\frac{\gamma\mu A}{D(CE)^{1/2}}\\
Ric(\hat{Y_4},\hat{Y_5})&=\frac{-\alpha\beta A}{B(DE)^{1/2}}+\frac{-\gamma\delta A}{C(DE)^{1/2}}.
\end{aligned}
\end{equation}
Then if we require $\alpha=\beta=\gamma=\delta=\eta=\mu=\rho=0$ the off-diagonal components also become 0. Thus, the metric will remain diagonal under the Ricci flow and we can restate $g$.\\
\\
\indent
Let $g(t)=A(t)\theta_1^2+B(t)\theta_2^2+C(t)\theta_3^2+D(t)\theta_4^2+E(t)\theta_5^2$ where $A(0)=\lambda_1,~B(0)=\lambda_2,~C(0)=\lambda_3,~D(0)=\lambda_4,$ and $E(0)=\lambda_5$. Then Ricci flow (\ref{Riccitensor}) reduces to the following system of equations:
\begin{align}\label{D5diffeqn}
\frac{dA}{dt}&=-\frac{A^2}{BC}\notag\\
\frac{dB}{dt}&=\frac{A}{C}\notag\\
\frac{dC}{dt}&=\frac{A}{B}\\
\frac{dD}{dt}&=0\notag\\
\frac{dE}{dt}&=4.\notag
\end{align}
Clearly, $D(t)=\lambda_4$ and $E(t)=4t+\lambda_5$.
To find the remaining equations we begin by calculating
\begin{equation*}\label{D5constants}
\frac{d}{dt}(AB)=\frac{d}{dt}(AC)=0.
\end{equation*}
Thus we have 
\begin{equation}\label{D5inverseconstants}
\frac{1}{B}=\frac{A}{\lambda_1\lambda_2}\text{ and }A=\frac{\lambda_1\lambda_3}{C}.
\end{equation}
So by (\ref{D5diffeqn})
\begin{align}
\frac{dC}{dt}&=\frac{A}{B}=\frac{A^2}{\lambda_1\lambda_2}=\frac{1}{\lambda_1\lambda_2}\Bigg(\frac{\lambda_1^2\lambda_3^2}{C^2}\Bigg),\label{D5dC}\\
C^2\frac{dC}{dt}&=\frac{\lambda_1\lambda_3^2}{\lambda_2}.\label{D5C^2}
\end{align}
Which leads to a solution for $C$,
\begin{align}
C^3&=\lambda_3^3+\frac{3\lambda_1\lambda_3^2}{\lambda_2}t\notag\\
C&=\lambda_3\Bigg(1+\frac{3\lambda_1}{\lambda_2\lambda_3}t\Bigg)^{1/3}.\label{D5C}
\end{align}
Thus by (\ref{D5inverseconstants}), we have the following solution to (\ref{D5diffeqn}):
\begin{equation}\label{D5solution}
\begin{aligned}
A&=\lambda_1\left(1+\frac{3\lambda_1}{\lambda_2\lambda_3}t\right)^{-1/3}\\
B&=\lambda_2\left(1+\frac{3\lambda_1}{\lambda_2\lambda_3}t\right)^{-1/3}\\
C&=\lambda_3\left(1+\frac{3\lambda_1}{\lambda_2\lambda_3}t\right)^{-1/3}\\
D&=\lambda_4\\
E&=4t+\lambda_5.
\end{aligned}
\end{equation}
\subsection{D11.}
Here we may choose a basis for the Lie Algebra $\{X_1,X_2,X_3,X_4,X_5\}$ such that the Lie bracket is of the form
\begin{align*}
[X_1,X_2]&=0 & [X_1,X_3]&=0 & [X_1,X_4]&=0 & [X_1,X_5]&=0 & [X_2,X_3]&=X_1\\
[X_2,X_4]&=0 & [X_2,X_5]&=X_3 & [X_3,X_4]&=0 & [X_3,X_5]&=-X_2 & [X_4,X_5]&=\epsilon X_1
\end{align*}
where $\epsilon=\pm1$.
Here we diagonalize the metric by letting $Y_i=\Lambda_i^kX_k$ with
\begin{equation*}
\Lambda=\left[\begin{array}{ccccc}
1&0&0&0&0\\a_1&1&0&0&0\\a_2&a_5&1&0&0\\a_3&a_6&a_8&1&0\\a_4&a_7&a_9&a_{10}&1
\end{array}\right].
\end{equation*}
We find that $\{Y_i\}$ satisfies the following bracket relations:
\begin{align*}
[Y_1,Y_2]&=0 & [Y_1,Y_3]&=0 & [Y_1,Y_4]&=0\\
[Y_1,Y_5]&=0 & [Y_2,Y_3]&=Y_1 & [Y_2,Y_4]&=\alpha Y_1\\
[Y_2,Y_5]&=\beta Y_1+\gamma Y_2+Y_3 & [Y_3,Y_4]&=\delta Y_1 & [Y_3,Y_5]&=\eta Y_1+(-1-\gamma^2)Y_2-\gamma Y_3\\
[Y_4,Y_5]&=\kappa Y_1+\rho Y_2+\sigma Y_3
\end{align*}
where $\alpha=a_{10}, \beta=a_1a_5-a_2+a_9, \gamma=-a_5, \delta=a_5a_{10}-a_6, \eta=a_1+a_1a_5^2-a_2a_5+a_5a_9-a_7, \kappa=a_1a_{10}+a_1a_5a_6-a_2a_6+a_6a_9-a_7a_{10}+\epsilon,\text{ and }\rho=-a_5a_6-a_{10}.$\\
\\
\indent
Define $\ds{\hat{Y_1}=\frac{Y_1}{\sqrt{A}}}$ as in (\ref{Yhat}). Let $W=w_1\hat{Y_1}+w_2\hat{Y_2}+w_3\hat{Y_3}+w_4\hat{Y_4}+w_5\hat{Y_5}$. Then by (\ref{Ricci}) we get the equations of the non-zero components of the Ricci tensor:
\begin{equation}\label{D11Ricci}
\begin{aligned}
Ric(\hat{Y_1},\hat{Y_1})&=\frac{A}{2BC}+\frac{\alpha^2A}{2BD}+\frac{\beta^2A}{2BE}+\frac{\delta^2A}{2CD}+\frac{\eta^2A}{2CE}+\frac{\kappa^2A}{2DE}\\
Ric(\hat{Y_2},\hat{Y_2})&=\frac{(-1-\gamma^2)^2B}{2CE}+\frac{\rho^2B}{2DE}+\frac{-A}{2BC}+\frac{-\alpha^2A}{2BD}+\frac{-\beta^2A}{2BE}+\frac{-C}{2BE}\\
Ric(\hat{Y_3},\hat{Y_3})&=\frac{C}{2BE}+\frac{\sigma^2C}{2DE}+\frac{-A}{2BC}+\frac{-\delta^2A}{2CD}+\frac{-\eta^2A}{2CE}+\frac{-(-1-\gamma^2)^2B}{2CE}\\
Ric(\hat{Y_4},\hat{Y_4})&=\frac{-\alpha^2A}{2BD}+\frac{-\delta^2A}{2CD}+\frac{-\kappa^2A}{2DE}+\frac{-\rho^2B}{2DE}+\frac{-\sigma^2C}{2DE}\\
Ric(\hat{Y_5},\hat{Y_5})&=\frac{-\beta^2A}{2BE}+\frac{-\gamma^2B}{2BE}+\frac{-C}{2BE}+\frac{-\eta^2A}{2CE}+\frac{-(-1-\gamma^2)^2B}{2CE}\\
&+\frac{-\gamma^2C}{2CE}+\frac{-\kappa^2A}{2DE}+\frac{-\rho^2B}{2DE}+\frac{-\sigma^2C}{2DE}+\frac{-\gamma B}{2BE}+\frac{-\gamma C}{2CE}\\
Ric(\hat{Y_2},\hat{Y_3})&=\frac{-\alpha\delta A}{D(BC)^{1/2}}+\frac{-\beta\eta A}{E(BC)^{1/2}}+\frac{-\gamma(-1-\gamma^2)B}{E(BC)^{1/2}}+\frac{\gamma C}{E(BC)^{1/2}}\\
Ric(\hat{Y_2},\hat{Y_4})&=\frac{\delta A}{C(BD)^{1/2}}+\frac{-\beta\kappa A}{E(BD)^{1/2}}+\frac{-\gamma\rho B}{E(BD)^{1/2}}+\frac{-\sigma C}{E(BD)^{1/2}}\\
Ric(\hat{Y_2},\hat{Y_5})&=\frac{\eta A}{C(BE)^{1/2}}+\frac{\alpha\kappa A}{D(BE)^{1/2}}\\
Ric(\hat{Y_3},\hat{Y_4})&=\frac{-\alpha A}{B(CD)^{1/2}}+\frac{-\eta\kappa A}{E(CD)^{1/2}}+\frac{-\rho(-1-\gamma^2)B}{E(CD)^{1/2}}+\frac{\gamma\sigma C}{E(CD)^{1/2}}\\
Ric(\hat{Y_3},\hat{Y_5})&=\frac{-\beta A}{B(CE)^{1/2}}+\frac{\delta\kappa A}{D(CE)^{1/2}}\\
Ric(\hat{Y_4},\hat{Y_5})&=\frac{-\alpha\beta A}{B(DE)^{1/2}}+\frac{-\delta\eta A}{C(DE)^{1/2}}.\\
\end{aligned}
\end{equation}
Then if we require $a_5=a_6=a_{10}$, $a_1=a_7$, and $a_2=a_9$, then $\alpha=\beta=\gamma=\delta=\eta=\rho=0$ and $\kappa=\epsilon$. So $\kappa^2=\epsilon^2=(\pm1)^2=1$. Therefore, the off-diagonal components are also forced to 0. Thus, the metric will remain diagonal under the Ricci flow and we can restate $g$.\\
\\
\indent
Let $g(t)=A(t)\theta_1^2+B(t)\theta_2^2+C(t)\theta_3^2+D(t)\theta_4^2+E(t)\theta_5^2$ where $A(0)=\lambda_1,~B(0)=\lambda_2,~C(0)=\lambda_3,~D(0)=\lambda_4,$ and $E(0)=\lambda_5$. Then Ricci flow (\ref{Riccitensor}) reduces to the following system of equations:
\begin{align}\label{D11diffeqn}
\frac{dA}{dt}&=-\frac{A^2}{BC}+-\frac{\kappa^2A^2}{DE}\notag\\
\frac{dB}{dt}&=-\frac{B^2}{CE}+\frac{A}{C}+\frac{C}{E}\notag\\
\frac{dC}{dt}&=-\frac{C^2}{BE}+\frac{A}{B}+\frac{B}{E}\\
\frac{dD}{dt}&=\frac{\kappa^2A}{E}\notag\\
\frac{dE}{dt}&=\frac{C}{B}+\frac{B}{C}+\frac{\kappa^2A}{D}.\notag
\end{align}
We restricted $\kappa$ in order that $\kappa^2=\epsilon^2=(\pm1)^2=1$, so we can then restate (\ref{D11diffeqn}):
\begin{align}\label{D11diffeqn1}
\frac{dA}{dt}&=-\frac{A^2}{BC}+-\frac{A^2}{DE}\notag\\
\frac{dB}{dt}&=-\frac{B^2}{CE}+\frac{A}{C}+\frac{C}{E}\notag\\
\frac{dC}{dt}&=-\frac{C^2}{BE}+\frac{A}{B}+\frac{B}{E}\\
\frac{dD}{dt}&=\frac{A}{E}\notag\\
\frac{dE}{dt}&=\frac{C}{B}+\frac{B}{C}+\frac{A}{D}.\notag
\end{align}
We calculate
\begin{align}
\frac{d}{dt}(A^2BCD^2)&=\frac{d}{dt}A^2E^2(B^2-C^2)=0\label{D11constant1}\\
\frac{d}{dt}(B-C)&=(B-C)(AE-B^2-2BC-C^2).\label{D11constant2}
\end{align}
Therefore, we may consider that if $B=C$ then $\ds{\frac{d}{dt}(B-C)=0}$. In contrast, if $B>C$ or $B<C$ at some time $t$, then it will remain so.\\
\\
\indent
Thus, by (\ref{D11constant2}) we can consider two cases. First, we let $\lambda_2=\lambda_3.$ Second, without a loss of generality, we let $\lambda_2>\lambda_3$ by the symmetry of $B$ and $C$.\\
\\
\textbf{Case 1:} Let $\lambda_2=\lambda_3$\\
\\
This condition gives the constants:
\begin{align}\label{D11.1constants}
B&=C\notag\\
A^2E^2(B^2-C^2)=\lambda_1^2\lambda_5^2(&\lambda_2^2-\lambda_3^2)=\lambda_1^2\lambda_5^2(0)=0\\
A^2BCD^2=\lambda_1^2\lambda_2&\lambda_3\lambda_4^2=\lambda_1^2\lambda_2^2\lambda_4^2.\notag
\end{align}
Then we can restate (\ref{D11diffeqn1}):
\begin{align}\label{D11.1diffeqn2}
\frac{dA}{dt}&=-\frac{A^2}{B^2}+-\frac{A^2}{DE}\notag\\
\frac{dB}{dt}&=\frac{A}{B}\notag\\
\frac{dD}{dt}&=\frac{A}{E}\\
\frac{dE}{dt}&=2+\frac{A}{D}\notag\\
B&=C.\notag
\end{align}
As such we know that
\begin{equation}\label{D11.1eqn}
A\leq\lambda_1,\quad B=C\geq\lambda_2,\quad D\geq\lambda_4,\quad E\geq\lambda_5+2t.
\end{equation}
We can then solve for A using (\ref{D11.1constants}) and (\ref{D11.1diffeqn2}):
\begin{align}\label{D11.1A}
\frac{dA}{dt}&=-\frac{A^2}{BC}+-\frac{A^2}{DE}\notag\\
&=-\frac{A^2}{B^2}+-\frac{A^2}{DE}\notag\\
&<-\frac{A^2}{B^2}\notag\\
&=-A^2\Bigg(\frac{A^2D^2}{\lambda_1^2\lambda_2^2\lambda_4^2}\Bigg)\notag\\
-\frac{1}{A^4}\frac{dA}{dt}&>\frac{1}{\lambda_1^2\lambda_2^2}\notag\\
\frac{1}{3}A^{-3}&>\frac{1}{\lambda_1^2\lambda_2^2}t+\frac{1}{3\lambda_1^3}\notag\\
A&<\lambda_1\Bigg(1+\frac{3\lambda_1}{\lambda_2^2}t\Bigg)^{-1/3}.
\end{align}
This solution for $A$ and \ref{D11.1eqn} allow us to solve for $D$:
\begin{align}\label{D11.1D}
\frac{dD}{dt}&=\frac{A}{E}\notag\\
&\leq\frac{\lambda_1\Bigg(1+\frac{3\lambda_1}{\lambda_2^2}t\Bigg)^{-1/3}}{2t+\lambda_5}\notag\\
&<\frac{\lambda_1}{2t+\lambda_5}\notag\\
\frac{dD}{dt}&<(\kappa_1t+\kappa_2)^{-4/3}\notag\\
D&<-3(\kappa_1t+\kappa_2)^{-1/3}+\lambda_4+\frac{3}{\kappa_2^{1/3}}\\
&<\lambda_4+\frac{3}{\kappa_2^{1/3}}\notag.
\end{align}
Thus, by the Monotone Convergence Theorem we observe that $D$ approaches a constant $K_1$.
Then, from (\ref{D11.1A}) follows a solution for B:
\begin{align}\label{D11.1B}
\frac{dB}{dt}&=\frac{A}{B}\notag\\
B\frac{dB}{dt}&<\lambda_1\Bigg(1+\frac{3\lambda_1}{\lambda_2^2}t\Bigg)^{-1/3}\notag\\
\frac{1}{2}B^2&<\frac{3}{2}\Bigg(\frac{\lambda_2^2}{3\lambda_1}\Bigg)\lambda_1\Bigg(1+\frac{3\lambda_1}{\lambda_2^2}t\Bigg)^{2/3}\notag\\
C=B&<\lambda_2\Bigg(1+\frac{3\lambda_1}{\lambda_2^2}t\Bigg)^{1/3}.
\end{align}
The solutions for $B$ and $D$ then allow us to find a value which $A$ is greater than:
\begin{align}\label{D11.1A1}
A^2&=\frac{\lambda_1^2\lambda_2^2\lambda_4^2}{B^2D^2}\notag\\
&>\frac{\lambda_1^2\lambda_2^2\lambda_4^2}{\lambda_2^2K_1^2}\Bigg(1+\frac{3\lambda_1}{\lambda_2^2}t\Bigg)^{-2/3}\notag\\
A&>\frac{\lambda_1\lambda_4}{K_1}\Bigg(1+\frac{3\lambda_1}{\lambda_2^2}t\Bigg)^{-1/3}.
\end{align}
We can solve for $C$ in a similar fashion:
\begin{align}\label{D11.1C}
\frac{dC}{dt}&=\frac{A}{C}\notag\\
C\frac{dC}{dt}&>\frac{\lambda_1\lambda_4}{K_1}\Bigg(1+\frac{3\lambda_1}{\lambda_2^2}t\Bigg)^{-1/3}\notag\\
\frac{1}{2}C^2&>\frac{3}{2}\Bigg(\frac{\lambda_2^2}{3\lambda_1}\Bigg)\frac{\lambda_1\lambda_4}{K_1}\Bigg(1+\frac{3\lambda_1}{\lambda_2^2}t\Bigg)^{2/3}+\eta\notag\\
C^2&>\frac{\lambda_2^2\lambda_4}{K_1}\Bigg(1+\frac{3\lambda_1}{\lambda_2^2}t\Bigg)^{2/3}+\eta\notag\\
\lambda_3^2&=\frac{\lambda_2^2\lambda_4}{K_1}+\eta\notag\\
\lambda_2^2=\lambda_3^2&=\frac{\lambda_2^2\lambda_4}{K_1}+\eta<\lambda_2^2+\eta\notag\\
\eta&>0\notag\\
B=C&>\Bigg(\frac{\lambda_4}{K_1}\Bigg)^{1/2}\lambda_2\Bigg(1+\frac{3\lambda_1}{\lambda_2^2}t\Bigg)^{1/3}.
\end{align}
Finally, we can also find $E$ from (\ref{D11.1diffeqn2}):
\begin{align}\label{D11.1E1}
\frac{dE}{dt}&=2+\frac{A}{D}\notag\\
&<2+\frac{\lambda_1}{\lambda_4}\notag\\
E&<\Bigg(2+\frac{\lambda_1}{\lambda_4}\Bigg)t+\lambda_5.
\end{align}
This gives the solution to (\ref{D11diffeqn1})
\begin{equation}\label{D11.1solution}
\begin{aligned}
\frac{\lambda_1\lambda_4}{K_1}\Bigg(1+\frac{3\lambda_1}{\lambda_2^2}t\Bigg)^{-1/3}<A(t)&<\lambda_1\Bigg(1+\frac{3\lambda_1}{\lambda_2^2}t\Bigg)^{-1/3}\\
\Bigg(\frac{\lambda_4}{K_1}\Bigg)^{1/2}\lambda_2\Bigg(1+\frac{3\lambda_1}{\lambda_2^2}t\Bigg)^{1/3}<C(t)=B(t)&<\lambda_2\Bigg(1+\frac{3\lambda_1}{\lambda_2^2}t\Bigg)^{1/3}\\
D(t)&\rightarrow K_1\\
2t+\lambda_5<E(t)&<\Bigg(2+\frac{\lambda_1}{\lambda_4}\Bigg)t+\lambda_5.
\end{aligned}
\end{equation}
This leads to:
\begin{equation}\label{D11.1solutionk}
\begin{aligned}
A(t)&\sim k_At^{-1/3}\\
C(t)=B(t)&\sim k_Bt^{1/3}\\
D(t)&\rightarrow K_1\\
E(t)&\sim k_Et^{-1/3}.
\end{aligned}
\end{equation}\\
\\
\textbf{Case 2:} Let $\lambda_2>\lambda_3$\\
\\
\begin{equation}\label{D11.2constants}
\begin{aligned}
A^2E^2(B^2-C^2)&=\lambda_1^2\lambda_5^2(\lambda_2^2-\lambda_3^2)\\
A^2BCD^2&=\lambda_1^2\lambda_2\lambda_3\lambda_4^2.
\end{aligned}
\end{equation}
By (\ref{D11diffeqn1}) we know that:
\begin{equation}\label{D11.2eqn}
A\leq\lambda_1,\quad D\geq\lambda_4,\quad E\geq\lambda_5+2t.
\end{equation}
Then,
\begin{align}
\frac{1}{2}\Bigg(\frac{1}{B+C}\frac{d(B+C)}{dt}+\frac{-1}{B-C}\frac{d(B-C)}{dt}\Bigg)&=\frac{1}{2}\Bigg(\frac{4}{E}\Bigg)\notag\\
&=\frac{2}{E}\notag\\
&=\frac{1}{E}\frac{dE}{dt}-\frac{1}{D}\frac{dD}{dt}.\notag\\
\frac{d}{dt}\Bigg(\frac{(B+C)D^2}{(B-C)E^{2}}\Bigg)&=0.\label{D11.2BCDE}
\intertext{Furthermore, since $\lambda_2>\lambda_3$,}
\frac{(B+C)D^2}{(B-C)E^2}=\frac{(\lambda_2+\lambda_3)\lambda_4^2}{(\lambda_2-\lambda_3)\lambda_5^2}
\end{align}
is positive.
Since $\ds{\frac{dE}{dt}>2+\frac{A}{D}}$ we know that $E>2t+\lambda_5$. Then $E^2>4t^2+4\lambda_5t+\lambda_5^2>4t^2$. Because $\lambda_2>\lambda_3$ and $D>\lambda_4$, we know that:
\begin{align}
\frac{(B+C)D^2}{(B-C)E^2}&=\frac{(\lambda_2+\lambda_3)\lambda_4^2}{(\lambda_2-\lambda_3)\lambda_5^2}\\
\frac{B+C}{B-C}&=\alpha\frac{E^2}{D^2}\\
&\geq\Bigg(\frac{(2t+\lambda_5)}{\frac{\lambda_1}{2}\ln(t+\frac{\lambda_5}{2})}\Bigg)^2\notag\\
&=\beta_1(t)\notag\\
&>1\\
(B+C)&>\beta_1(t)(B-C)\\
\frac{B}{C}&<\frac{\beta_1(t)+1}{\beta_1(t)-1}.
\intertext{Then we see that for any $\delta_1>0$ there is some $T$ such that for all $t>T$,}
1<\frac{B}{C}&<\delta_1.
\intertext{This leads us to observe that}
\frac{B}{C}&\rightarrow1.
\end{align}
This allows us to solve (\ref{D11diffeqn1}) the same way we did in Case 1. If we follow the same method we generate the following solution, where $k_1$ and $k_2$ are constants and $\ds{k_1=\frac{\lambda_2\lambda_3\lambda_4}{\delta_1K_1}}$ and $k_2=\delta_1\lambda_2\lambda_3$.
This leads to:
\begin{equation}\label{D11.2solutionk}
\begin{aligned}
A(t)&\sim k_At^{-1/3}\\
C(t)\lesssim B(t)&\sim k_Bt^{1/3}\\
D(t)&\rightarrow K_1\\
E(t)&\sim k_Et.
\end{aligned}
\end{equation}
\section{Conclusion}
This project demonstrates significant differences in the behavior of Ricci flow in five dimensions from what has been observed in two, three, and four dimensions. Much of the lower dimensional behavior remained linear or inverse cubic. However, while classes D1, D5, and D11 also demonstrated inverse cubic behavior, classes D2 and D3 actually developed solutions raised to inverse seventh and inverse eleventh powers. This is indicative of the changes Ricci flow goes through when applied to higher dimensions, and suggests that greater oddities will occur in even higher dimensions. It would be valuable to research other five-dimensional geometries and discover whether the cubic behavior is retained or not. Research into backwards Ricci flow in five-dimensions as well as looking into higher dimensions will also be able to provide clues as to the more general nature of Ricci flow across several different dimensions.


\begin{thebibliography}{}

\bibitem{Be} Bell, Thomas, \emph{Backward Ricci Flow of Compact Locally Homogeneous Geometries on $4$-Manifolds}. arXiv:1507.08988

\bibitem{B} A. Besse, \emph{Einstein manifolds}, Springer-Verlag, Berlin, 1987. p.183.

\bibitem{BD} Boyce, William E., DiPrima, Richard C., \emph{Elementary Differential Equations, 10th Edition}, John Wiley and Sons, Inc., Hoboken, New Jersey, 2012. p.77.

\bibitem{CSC} Cao, Xiaodong; Saloff-Coste, Laurent, \emph{Backward Ricci flow on locally homogeneous three-manifolds}. Comm. Anal. Geom., \textbf{17(2)} (2009), 305-325.

\bibitem{DF} Diatta, Andre; Foreman, Brendan, \emph{Lattices in contact Lie groups and 5-dimensional contact solvmanifolds}. arXiv:0904.3113

\bibitem{Ha} Hamilton, Richard, \emph{Three-manifolds with positive Ricci curvature}. J. Differential Geom. \textbf{17} (1982), 255–-306.

\bibitem{IJ}  Isenberg, James; Jackson, Martin, \emph{Ricci flow of locally homogeneous geometries on closed manifolds}. J. Differential Geom. \textbf{35} (1992), 733-741.

\bibitem{IJL} Isenberg, James; Jackson, Martin; Lu, Peng, \emph{Ricci flow on locally homogeneous closed $4$-manifolds}. Comm. Anal. Geom., \textbf{14(2)} (2006), 345-386.

\bibitem{Pe} Perelman, Grigori, \emph{The entropy formula for the Ricci flow and its geometric applications}. arXiv:math/0211159

\bibitem{Per} Perelman, Grigori, \emph{Ricci flow with surgery on three-manifolds}. arXiv:math/0303109

\bibitem{Pere} Perelman, Grigori, \emph{Finite extinction time for the solutions to the Ricci flow on certain three-manifolds}. 	arXiv:math/0307245

\end{thebibliography}
\end{document}